\pdfoutput=1 

\documentclass[hidelinks]{ifacconf}

\makeatletter
\let\old@ssect\@ssect 
\makeatother

\usepackage{graphicx}      
\usepackage{natbib}        

\usepackage{amsmath,amsfonts,amssymb}
\usepackage{mathtools}
\usepackage{lmodern}

\usepackage{hyperref}
\usepackage{xurl} 
        
\usepackage{algorithm}
\usepackage{algorithmicx}
\usepackage{algpseudocode}
\algnewcommand\algorithmicinput{\textbf{Input:}}
\algnewcommand\Input{\item[\algorithmicinput]}
\algnewcommand\algorithmicoutput{\textbf{Output:}}
\algnewcommand\Output{\item[\algorithmicoutput]}

\newcommand{\N}{\mathbb{N}}
\newcommand{\R}{\mathbb{R}}

\DeclareMathOperator{\dom}{dom}
\DeclareMathOperator{\im}{im}
\newcommand{\norm}[1]{\Vert #1 \Vert}
\newcommand{\ip}[2]{\langle #1, #2 \rangle}
\DeclareMathOperator*{\argmin}{argmin}

\newcommand{\setmapsto}{\multimap}

\newcommand{\proj}[2]{P_{#1}(#2)}

\newcommand{\gencone}[4]{{#1}_{#2}^{#4}(#3)} 
\newcommand{\tancone}[2]{\gencone{T}{#1}{#2}{}} 

\newcommand{\ushort}[1]{\mkern 0.9mu\underline{\mkern-0.9mu#1\mkern-0.9mu}\mkern 0.9mu}

\newcommand{\tp}{\top}
\DeclareMathOperator{\tr}{tr}
\DeclareMathOperator{\rank}{rank}
\DeclareMathOperator{\diag}{diag}
\newcommand{\st}{\mathrm{St}}

\newcommand{\ppgd}{\mathrm{P}^2\mathrm{GD}}

\makeatletter
\def\@ssect#1#2#3#4#5#6{%
  \NR@gettitle{#6}
  \old@ssect{#1}{#2}{#3}{#4}{#5}{#6}
}
\makeatother

\edef\endfrontmatter{%
  \unexpanded\expandafter{\endfrontmatter}
  \noexpand\endNoHyper 
}
\begin{document}
\begin{frontmatter}

\title{Retractions on closed sets\thanksref{footnoteinfo}} 

\thanks[footnoteinfo]{This work was supported by the ERC grant \#786854 G-Statistics from the European Research Council under the European Union's Horizon 2020 research and innovation program and by the French government through the 3IA C\^{o}te d'Azur Investments ANR-19-P3IA-0002 managed by the National Research Agency.}

\author[First]{Guillaume Olikier} 

\address[First]{Universit\'{e} C\^{o}te d'Azur and Inria, Epione Project Team, 2004 route des Lucioles - BP 93, 06902 Sophia Antipolis Cedex, France (\href{mailto:guillaume.olikier@inria.fr}{\nolinkurl{guillaume.olikier@inria.fr}})}

\begin{abstract}
On a manifold or a closed subset of a Euclidean vector space, a retraction enables to move in the direction of a tangent vector while staying on the set. Retractions are a versatile tool to perform computational tasks such as optimization, interpolation, and numerical integration. This paper studies two known definitions of retraction on a closed subset of a Euclidean vector space, one being weaker than the other. Specifically, it shows that, in the context of constrained optimization, the weaker definition should be preferred as it inherits the main property of the other while being less restrictive.
\end{abstract}


\end{frontmatter}

\section{Introduction}
\label{sec:Introduction}
In topology, a retraction is a continuous map from a topological space into a subspace whose restriction to this subspace is the identity map on the subspace \cite[\S 1.2]{tomDieck}.
The retractions on manifolds or closed subsets of Euclidean vector spaces considered in this paper are first-order retractions: they are retractions in the sense of topology that satisfy a first-order condition (see Section~\ref{sec:ProjectiveRetractionsClosedSets}). These retractions enable to generate trajectories on the set on which they are defined with prescribed initial position and velocity.

A retraction on a manifold is a smooth map from the tangent bundle, i.e., the union of all tangent spaces, onto the manifold that satisfies some conditions; see \citet[\S 4.1]{AbsilMahonySepulchre} or \citet[\S\S 3.6 and 8.7]{Boumal} for a precise definition and \citet[\S 4.10]{AbsilMahonySepulchre} for historical notes.
Enabling to move from a point on the manifold in the direction of a tangent vector while staying on the manifold, retractions offer an efficient alternative to the Riemannian exponential map to perform computational tasks on manifolds, such as optimization \citep{AbsilMalick,AbsilOseledets2015}, interpolation \citep{SeguinKressner,ZimmermannBergmann}, and numerical integration \citep{SeguinCerutiKressner}.

To enlarge their scope of applicability, the concept of retraction has been extended from manifolds to arbitrary closed subsets of Euclidean vector spaces, which merely requires the replacement of the tangent space with the tangent cone. In the context of optimization on the determinantal variety \cite[Lecture~9]{Harris}
\begin{equation*}
\R_{\le r}^{m \times n} \coloneq \{X \in \R^{m \times n} \mid \rank X \le r\},
\end{equation*}
$m$, $n$, and $r$ being positive integers such that $r < \min\{m, n\}$, two such extensions are proposed in \citet[Definition~2.4]{SchneiderUschmajew2015} and \citet[\S 3.1.2]{HosseiniUschmajew2019}. Contrary to the latter, the former requires every trajectory generated by the retraction to be continuous, thereby defining a stronger concept.

On a submanifold of a Euclidean vector space, the projection map defines a retraction called the projective retraction \cite[Proposition~5]{AbsilMalick}. On a closed subset of a Euclidean vector space, the projection map is suggested in \citet[\S 2.4.1]{SchneiderUschmajew2015} to define a retraction in the sense of \citet[\S 3.1.2]{HosseiniUschmajew2019} if the set satisfies a regularity condition, called geometric derivability, satisfied by all algebraic varieties. Furthermore, retractions in the sense of \citet[Definition~2.4]{SchneiderUschmajew2015} ensure a property playing an instrumental role in the design and analysis of algorithms for constrained optimization \cite[Proposition~2.8]{SchneiderUschmajew2015}.

This paper aims at showing that, in the context of constrained optimization, the definition of retraction proposed in \citet[\S 3.1.2]{HosseiniUschmajew2019} should be preferred to that proposed in \citet[Definition~2.4]{SchneiderUschmajew2015} as it inherits the main property while being less restrictive. Specifically, three contributions are brought to the study of retractions on closed subsets of Euclidean vector spaces. In Section~\ref{sec:ProjectiveRetractionsClosedSets}, the projection map is proven to define a retraction in the sense of \citet[\S 3.1.2]{HosseiniUschmajew2019} on geometrically derivable sets. In Section~\ref{sec:ProjectiveRetractionsDeterminantalVariety}, it is shown that the projection onto the determinantal variety $\R_{\le r}^{m \times n}$ does not define a retraction in the sense of \citet[Definition~2.4]{SchneiderUschmajew2015}. In Section~\ref{sec:RetractedLineSearchClosedSet}, the property given in \citet[Proposition~2.8]{SchneiderUschmajew2015} is proven to be satisfied for retractions in the sense of \citet[\S 3.1.2]{HosseiniUschmajew2019}. These results are important for the convergence analysis of algorithms aiming at minimizing a differentiable function on a nonempty closed subset of a Euclidean vector space, such as those proposed in \citet[Algorithm~3]{SchneiderUschmajew2015} and \citet[Algorithm~2]{GaoPengYuan} for low-rank matrix and tensor varieties. As a matter of fact, the retractions on low-rank tensor varieties proposed in \citet[\S 5]{Kutschan2018} and \citet[Proposition~2]{GaoPengYuan} satisfy the definition from \citet[\S 3.1.2]{HosseiniUschmajew2019} but are not proven to satisfy that from \citet[Definition~2.4]{SchneiderUschmajew2015}.

The necessary background material in variational analysis is reviewed in Section~\ref{sec:ElementsVariationalAnalysis}. Conclusions are drawn in Section~\ref{sec:Conclusion}.
Throughout the paper, $\mathcal{E}$ is a Euclidean vector space with inner product $\ip{\cdot}{\cdot\cdot}$ and induced norm $\norm{\cdot}$, and $C$ is a nonempty closed subset of $\mathcal{E}$.

\section{Elements of variational analysis}
\label{sec:ElementsVariationalAnalysis}
This section, mostly based on \citet{RockafellarWets}, reviews background material in variational analysis that is used in the rest of the paper. Section~\ref{subsec:OuterSemicontinuityLocalBoundednessCorrespondence} is about outer semicontinuity and local boundedness of a correspondence. Section~\ref{subsec:ProjectionMap} concerns the projection map and its main properties. Section~\ref{subsec:TangentCone} reviews the concept of tangent cone on which the concept of retraction relies.
A nonempty subset $K$ of $\mathcal{E}$ is called a \emph{cone} if $x \in K$ implies $\alpha x \in K$ for all $\alpha \in [0, \infty)$ \cite[\S 3B]{RockafellarWets}.

\subsection{Outer semicontinuity and local boundedness of a correspondence}
\label{subsec:OuterSemicontinuityLocalBoundednessCorrespondence}
A \emph{correspondence}, or a \emph{set-valued map}, is a triplet $F \coloneq (A_1, A_2, G)$ where $A_1$ and $A_2$ are sets respectively called the \emph{set of departure} and the \emph{set of destination} of $F$, and $G$ is a subset of $A_1 \times A_2$ called the \emph{graph} of $F$. If $F \coloneq (A_1, A_2, G)$ is a correspondence, written $F : A_1 \setmapsto A_2$, then the \emph{image} of $x \in A_1$ by $F$ is $F(x) \coloneq \{y \in A_2 \mid (x, y) \in G\}$ and the \emph{domain} of $F$ is $\dom F \coloneq \{x \in A_1 \mid F(x) \ne \emptyset\}$. A \emph{selection} of a correspondence $F : A_1 \setmapsto A_2$ is a function $f : \dom F \to A_2$ such that $f(x) \in F(x)$ for all $x \in \dom F$.

Let $\mathcal{E}_1$ and $\mathcal{E}_2$ be two Euclidean vector spaces. A correspondence $F : \mathcal{E}_1 \setmapsto \mathcal{E}_2$ is said to be \emph{outer semicontinuous} at $x \in \mathcal{E}_1$ if, for every sequence $(x_i)_{i \in \N}$ in $\mathcal{E}_1$ converging to $x$ and every sequence $(y_i)_{i \in \N}$ in $\mathcal{E}_2$ converging to $y$ such that, for all $i \in \N$, $y_i \in F(x_i)$, it holds that $y \in F(x)$ \cite[Definition~5.4]{RockafellarWets}. Furthermore, $F$ is said to be \emph{locally bounded} at $x \in \mathcal{E}_1$ if there exists $\rho \in (0, \infty)$ such that $F(B_{\mathcal{E}_1}[x, \rho])$ is bounded, $B_{\mathcal{E}_1}[x, \rho]$ being the closed ball of center $x$ and radius $\rho$ in $\mathcal{E}_1$, and $F$ is said to be locally bounded if it is locally bounded at every point of $\mathcal{E}_1$ \cite[Definition~5.14]{RockafellarWets}. Let $\mathcal{E}_3$ be a Euclidean vector space and $H : \mathcal{E}_2 \setmapsto \mathcal{E}_3$ be a correspondence. If $F$ and $H$ are outer semicontinuous and $F$ is locally bounded, then $H \circ F$ is outer semicontinuous \cite[Proposition~5.52(b)]{RockafellarWets}.

\subsection{Projection map}
\label{subsec:ProjectionMap}
The distance from $x \in \mathcal{E}$ to $C$ is $d(x, C) \coloneq \min_{y \in C} \norm{x-y}$ and the projection of $x$ onto $C$ is $\proj{C}{x} \coloneq \argmin_{y \in C} \norm{x-y}$ which is a nonempty compact set \cite[Example~1.20]{RockafellarWets}. When $\proj{C}{x}$ is a singleton, it is identified with its unique element. The projection map $P_C$ is outer semicontinuous and locally bounded \cite[Example~5.23(a)]{RockafellarWets}.

\subsection{Tangent cone}
\label{subsec:TangentCone}
A vector $v \in \mathcal{E}$ is said to be \emph{tangent} to $C$ at $x \in C$ if there exist sequences $(x_i)_{i \in \N}$ in $C$ converging to $x$ and $(t_i)_{i \in \N}$ in $(0, \infty)$ such that the sequence $(\frac{x_i-x}{t_i})_{i \in \N}$ converges to~$v$ \cite[Definition~6.1]{RockafellarWets}. The set of all tangent vectors to $C$ at $x \in C$ is a closed cone \cite[Proposition~6.2]{RockafellarWets} called the \emph{tangent cone} to $C$ at $x$ and denoted by $\tancone{C}{x}$.
If $C$ is a submanifold of $\mathcal{E}$ locally around $x \in C$, then $\tancone{C}{x}$ coincides with the tangent space to $C$ at $x$ \cite[Example~6.8]{RockafellarWets}. Given $x \in C$, a vector $v \in \tancone{C}{x}$ is said to be \emph{derivable} if there exist $\varepsilon \in (0, \infty)$ and $\gamma : [0, \varepsilon] \to C$ such that $\gamma(0) = x$ and $\gamma'(0) = v$. The set $C$ is said to be \emph{geometrically derivable} at $x \in C$ if every $v \in \tancone{C}{x}$ is derivable. For example, every real algebraic variety is geometrically derivable at each of its points \cite[Proposition~2]{OSheaWilson2004}.

As pointed out in \citet[\S 6A]{RockafellarWets}, for every $x \in C$, the vectors $v \in \tancone{C}{x}$ can be characterized based on the derivability at zero of the function
\begin{equation*}
[0, \infty) \to \R : t \mapsto d(x+tv, C).
\end{equation*}
Specifically, for every $x \in C$,
\begin{align*}
\tancone{C}{x} = \left\{v \in \mathcal{E} \mid \liminf_{t \to 0^+} \frac{d(x+tv, C)}{t} = 0\right\}
\end{align*}
and the subset of $\tancone{C}{x}$ containing the derivable tangent vectors is
\begin{align*}
\left\{v \in \mathcal{E} \mid \lim_{t \to 0^+} \frac{d(x+tv, C)}{t} = 0\right\}.
\end{align*}

In this paper, $\R^{m \times n}$ is endowed with the Frobenius inner product $(X, Y) \mapsto \tr Y^\tp X$. The tangent cone to $\R_{\le r}^{m \times n}$ at $X \in \R_{\le r}^{m \times n}$ can be described explicitly based on orthonormal bases of $\im X$, $\im X^\tp$, and their orthogonal complements \cite[Theorem~3.2]{SchneiderUschmajew2015}. 
For every $p, q \in \N$, $0_{p \times q}$ is the zero matrix in $\R^{p \times q}$, $0_p \coloneq 0_{p \times p}$, $I_p$ is the identity matrix in $\R^{p \times p}$, if $q \ge p$, then
\begin{equation*}
\st(p, q) \coloneq \{U \in \R^{q \times p} \mid U^\tp U = I_p\}
\end{equation*}
is a Stiefel manifold, and $\mathrm{O}(p) \coloneq \st(p, p)$ is an orthogonal group \cite[\S 3.3.2]{AbsilMahonySepulchre}. Let $\ushort{r} \coloneq \rank X$, $U \in \st(\ushort{r}, m)$, $U_\perp \in \st(m-\ushort{r}, m)$, $V \in \st(\ushort{r}, n)$, $V_\perp \in \st(n-\ushort{r}, n)$, $\im U = \im X$, $\im U_\perp = (\im X)^\perp$, $\im V = \im X^\tp$, and $\im V_\perp = (\im X^\tp)^\perp$.
Then,
\begin{equation}
\label{eq:TanConeDeterminantalVariety}
\tancone{\R_{\le r}^{m \times n}}{X} = [U \; U_\perp] \begin{bmatrix} \R^{\ushort{r} \times \ushort{r}} & \R^{\ushort{r} \times n-\ushort{r}} \\ \R^{m-\ushort{r} \times \ushort{r}} & \R_{\le r-\ushort{r}}^{m-\ushort{r} \times n-\ushort{r}} \end{bmatrix} [V \; V_\perp]^\tp.
\end{equation}

\section{Projective retractions on closed sets}
\label{sec:ProjectiveRetractionsClosedSets}
As pointed out in Section~\ref{sec:Introduction}, for every submanifold $\mathcal{M}$ of~$\mathcal{E}$, the projection onto $\mathcal{M}$ defines a retraction on $\mathcal{M}$ called the \emph{projective retraction} \cite[Proposition~5]{AbsilMalick}. In this section, the projection onto $C$ is proven to define a retraction on $C$ in the sense of \citet[\S 3.1.2]{HosseiniUschmajew2019} if $C$ is geometrically derivable, as suggested in \citet[\S 2.4.1]{SchneiderUschmajew2015}.

\begin{defn}[{\citet[\S 3.1.2]{HosseiniUschmajew2019}}]
\label{defi:Retraction}
A \emph{retraction} on $C$ is a function
\begin{equation*}
R : \bigcup_{x \in C} \{x\} \times \tancone{C}{x} \to C
\end{equation*}
such that, for all $x \in C$ and $v \in \tancone{C}{x}$,
\begin{equation*}
\lim_{t \to 0^+} \frac{R(x ,tv)-(x+tv)}{t} = 0.
\end{equation*}
\end{defn}

Thus, if $R$ is a retraction on $C$, then for every $x \in C$ and $v \in \tancone{C}{x}$ the function
\begin{equation}
\label{eq:CurveFromRetraction}
\hat{R} : [0, \infty) \to C : t \mapsto R(x, tv)
\end{equation}
satisfies
\begin{align}
\label{eq:Retraction}
\hat{R}(0) = x,&&
\hat{R}'(0) = v.
\end{align}
In contrast, a retraction $R$ in the sense of \citet[Definition~2.4]{SchneiderUschmajew2015} must further ensure that, for every $x \in C$ and $v \in \tancone{C}{x}$, the function $\hat{R}$ in~\eqref{eq:CurveFromRetraction} is continuous.

As announced in Section~\ref{sec:Introduction}, a retraction $R$ in the sense of Definition~\ref{defi:Retraction} is a retraction in the sense of topology. Indeed, $C$ can be identified with $C \times \{0\}$ by identifying every $x \in C$ with $(x, 0)$. Furthermore, since $0 \in \tancone{C}{x}$ for all $x \in C$, the set $C \times \{0\}$ is a subset of $\bigcup_{x \in C} \{x\} \times \tancone{C}{x}$. Therefore, the first equation in~\eqref{eq:Retraction} ensures that the restriction of $R$ to $C \times \{0\}$ is the identity map on $C \times \{0\}$. The first-order condition mentioned in Section~\ref{sec:Introduction} is the second equation in~\eqref{eq:Retraction}.

\begin{prop}
\label{prop:ProjectiveRetractionGeometricDerivability}
If $C$ is geometrically derivable at each of its points, then every function
\begin{equation*}
R : \bigcup_{x \in C} \{x\} \times \tancone{C}{x} \to C
\end{equation*}
such that, for all $x \in C$ and $v \in \tancone{C}{x}$,
\begin{equation*}
R(x, v) \in \proj{C}{x+v}
\end{equation*}
is a retraction on $C$ called a \emph{projective retraction} on $C$.
\end{prop}

\begin{pf}
Assume that $C$ is geometrically derivable at each of its points. Let $x \in C$ and $v \in \tancone{C}{x}$. Then, there exist $\varepsilon \in (0, \infty)$ and $\gamma : [0, \varepsilon] \to C$ such that $\gamma(0) = x$ and $\gamma'(0) = v$. Thus, for all $t \in (0, \varepsilon]$,
\begin{align*}
\frac{\norm{R(x, tv)-(x+tv)}}{t}
&= \frac{d(x+tv, C)}{t}\\
&\le \frac{\norm{\gamma(t)-(x+tv)}}{t}\\
&= \frac{\norm{\gamma(t)-\gamma(0)-t\gamma'(0)}}{t},
\end{align*}
and the result follows by letting $t$ tend to zero.
\hfill$\qed$
\end{pf}

In general, a projective retraction cannot be expected to be a retraction in the sense of \citet[Definition~2.4]{SchneiderUschmajew2015}. Indeed, for a projective retraction on $C$ to be a retraction in the sense of \citet[Definition~2.4]{SchneiderUschmajew2015}, it would be necessary that, for every $x \in C$ and $v \in \tancone{C}{x}$, a continuous selection of the correspondence
\begin{equation}
\label{eq:ProjectionTangent}
[0, \infty) \setmapsto C : t \mapsto \proj{C}{x+tv}
\end{equation}
exists, which is not the case in general. Indeed, in view of the properties reviewed in Sections~\ref{subsec:OuterSemicontinuityLocalBoundednessCorrespondence} and~\ref{subsec:ProjectionMap}, the correspondence defined in~\eqref{eq:ProjectionTangent} is outer semicontinuous. However, in general, outer semicontinuous correspondences cannot be expected to admit continuous selections \cite[p.~188]{RockafellarWets}. For $C = \R_{\le r}^{m \times n}$, it is proven in Section~\ref{sec:ProjectiveRetractionsDeterminantalVariety} that, for every $x \in C \setminus \{0\}$ and $t_* \in (0, \infty)$, there exists $v \in \tancone{C}{x}$ such that every selection of the correspondence defined in~\eqref{eq:ProjectionTangent} is discontinuous at $t_*$.

\section{Projective retractions on the determinantal variety}
\label{sec:ProjectiveRetractionsDeterminantalVariety}
\begin{figure}[h!]
\begin{center}
\vspace*{-2mm}
\includegraphics[width=\linewidth]{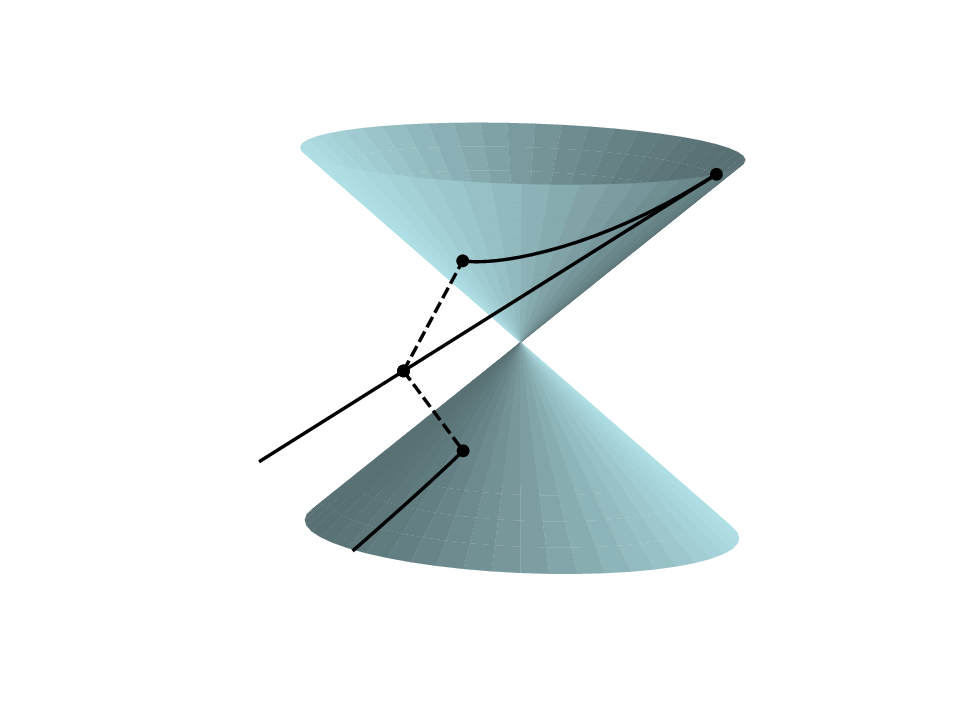}
\vspace*{-11mm}
\end{center}
\caption{Illustration of the proof of Proposition~\ref{prop:ProjectiveRetractionDeterminantalVariety}.
This proof can be adapted to show that Proposition~\ref{prop:ProjectiveRetractionDeterminantalVariety} still holds if $\R^{m \times n}$ and $\R_{\le r}^{m \times n}$ are replaced with $\mathrm{S}(n) \coloneq \{X \in \R^{n \times n} \mid X^\tp = X\}$ and $\mathrm{S}_{\le r}(n) \coloneq \{X \in \mathrm{S}(n) \mid \rank X \le r\}$, respectively. Moreover, this result can be illustrated for $n = 2$ and $r = 1$ by taking advantage of two facts. First, the map $\varphi : \R^3 \to \mathrm{S}(2) : x \mapsto \frac{1}{\sqrt{2}} \left[\begin{smallmatrix} x_3-x_1 & x_2 \\ x_2 & x_3+x_1 \end{smallmatrix}\right]$ is a bijection. Second, $\mathrm{S}_{\le 1}(2)$ is the image of the cone $C \coloneq \{x \in \R^3 \mid x_3^2 = x_1^2+x_2^2\}$ under $\varphi$.
Given $x \in C \setminus \{0\}$, let $Z$ be the tangent vector to $\mathrm{S}_{\le 1}(2)$ at $\varphi(x)$ defined in the proof of Proposition~\ref{prop:ProjectiveRetractionDeterminantalVariety}. Then, $Z \in \mathrm{S}(2)$. Define $z \coloneq \varphi^{-1}(Z)$.
The point $x$, the half-line $\{x+tz \mid t \in [0, \infty)\}$, and $\bigcup_{t \in [0, \infty)} \proj{C}{x+tz}$ are represented in the figure. The set $\proj{C}{x+tz}$ is a singleton for every $t \in [0, \infty) \setminus \{t_*\}$, but $\proj{C}{x+t_*z}$ contains two elements, one on the upper part of the cone, and the other on the lower part of the cone. The segments joining $x+t_*z$ to each of its two projections onto $C$ are represented in dashed line. This figure illustrates that every selection of $[0, \infty) \setmapsto C : t \mapsto \proj{C}{x+tz}$ is discontinuous at $t_*$.}
\label{fig:ProjectiveRetractionDoubleCone}
\end{figure}

This section concerns projective retractions on the determinantal variety $\R_{\le r}^{m \times n}$ which exist by Proposition~\ref{prop:ProjectiveRetractionGeometricDerivability}. Based on the description of the tangent cone to $\R_{\le r}^{m \times n}$ reviewed in~\eqref{eq:TanConeDeterminantalVariety}, Proposition~\ref{prop:ProjectiveRetractionDeterminantalVariety} implies that a projective retraction on $\R_{\le r}^{m \times n}$ is not a retraction in the sense of \citet[Definition~2.4]{SchneiderUschmajew2015}. The proof is illustrated in Figure~\ref{fig:ProjectiveRetractionDoubleCone}.

\begin{prop}
\label{prop:ProjectiveRetractionDeterminantalVariety}
For every $X \in \R_{\le r}^{m \times n} \setminus \{0_{m \times n}\}$ and $t_* \in (0, \infty)$, there exists $Z \in \tancone{\R_{\le r}^{m \times n}}{X}$ such that every selection of the correspondence
\begin{equation}
\label{eq:ProjectionTangentDeterminantalVariety}
[0, \infty) \setmapsto \R_{\le r}^{m \times n} : t \mapsto \proj{\R_{\le r}^{m \times n}}{X+tZ}
\end{equation}
is discontinuous at $t_*$.
\end{prop}

\begin{pf}
Consider an SVD of $X \in \R_{\le r}^{m \times n}$:
\begin{equation*}
X = U \diag(\sigma_1, \dots, \sigma_{\ushort{r}}, 0_{r-\ushort{r}+1}) V^\tp,
\end{equation*}
where $\ushort{r} \coloneq \rank X$, $\sigma_1 \ge \dots \ge \sigma_{\ushort{r}} > 0$ are the nonzero singular values of $X$, $U \in \st(r+1, m)$, and $V \in \st(r+1, n)$. (The $\diag$ operator returns a block diagonal matrix \cite[\S 1.3.1]{GolubVanLoan} whose diagonal blocks are the arguments given to the operator.) Thus, for every $i \in \{1, \dots, \ushort{r}\}$, the columns $U_{:, i}$ and $V_{:, i}$ are respectively left and right singular vectors associated with the singular value $\sigma_i$.
Let $t_* \in (0, \infty)$ and define
\begin{equation*}
Z \coloneq
\frac{\sigma_{\ushort{r}}}{t_*}
U
\diag\left(0_{\ushort{r}-1}, \left[\begin{smallmatrix} -1 & 1/2 \\ 1/2 & 0 \end{smallmatrix}\right], {\textstyle \frac{3}{4}} I_{r-\ushort{r}}\right)
V^\tp
\in \tancone{\R_{\le r}^{m \times n}}{X}.
\end{equation*}
Let $t \in (0, \infty)$ and define $\tau \coloneq t/t_*$. Then,
\begin{align*}
&X+tZ =\\
&U
\diag\left(\sigma_1, \dots, \sigma_{\ushort{r}-1}, \sigma_{\ushort{r}} \left[\begin{smallmatrix} 1-\tau & \tau/2 \\ \tau/2 & 0 \end{smallmatrix}\right], {\textstyle \frac{3}{4}} \tau \sigma_{\ushort{r}} I_{r-\ushort{r}}\right)
V^\tp
\end{align*}
and the nonzero singular values of $X+tZ$ are $\sigma_1, \dots, \sigma_{\ushort{r}-1}$, $\frac{3}{4} \tau \sigma_{\ushort{r}}$, and the absolute values of the eigenvalues of $\left[\begin{smallmatrix} 1-\tau & \tau/2 \\ \tau/2 & 0 \end{smallmatrix}\right]$, i.e.,
\begin{equation*}
\lambda_\pm(\tau) \coloneq \frac{1-\tau\pm\sqrt{(1-\tau)^2+\tau^2}}{2},
\end{equation*}
multiplied by $\sigma_{\ushort{r}}$. Thus, $X+tZ$ has rank $r+1$ and, by the Eckart--Young theorem \citep{EckartYoung1936}, projecting this matrix onto $\R_{\le r}^{m \times n}$ amounts to compute a truncated SVD of rank $r$.
With
\begin{equation*}
w_\pm(\tau) \coloneq \begin{bmatrix} \lambda_\pm(\tau) \\ \tau/2 \end{bmatrix} / \sqrt{\lambda_\pm(\tau)^2+\tau^2/4},
\end{equation*}
it holds that $[w_+(\tau) \; w_-(\tau)] \in \mathrm{O}(2)$ and the following is an eigendecomposition:
\begin{align*}
&\begin{bmatrix} 1-\tau & \tau/2 \\ \tau/2 & 0 \end{bmatrix} =\\
&\;[w_+(\tau) \; w_-(\tau)] \diag(\lambda_+(\tau), \lambda_-(\tau)) [w_+(\tau) \; w_-(\tau)]^\tp.
\end{align*}
Thus,
\begin{align*}
&X+tZ =\\
&\begin{bmatrix} U_{:, 1:\ushort{r}-1} & U_{:, \ushort{r}:\ushort{r}+1}[w_+(\tau) \; w_-(\tau)] & U_{:, \ushort{r}+2:r+1} \end{bmatrix}\\
&\,\diag(\sigma_1, \dots, \sigma_{\ushort{r}-1}, \sigma_{\ushort{r}} \lambda_+(\tau), \sigma_{\ushort{r}} \lambda_-(\tau), {\textstyle \frac{3}{4}} \tau \sigma_{\ushort{r}} I_{r-\ushort{r}})\\
&\begin{bmatrix} V_{:, 1:\ushort{r}-1} & V_{:, \ushort{r}:\ushort{r}+1}[w_+(\tau) \; w_-(\tau)] & V_{:, \ushort{r}+2:r+1} \end{bmatrix}^\tp.
\end{align*}
If $\tau \in (\frac{12}{17}, 1)$, then $1 > \frac{3}{4} \tau > \lambda_+(\tau) > -\lambda_-(\tau) > 0$. If $\tau \in (1, \frac{4}{3})$, then $1 > \frac{3}{4} \tau > -\lambda_-(\tau) > \lambda_+(\tau) > 0$.
Hence, a reduced SVD \cite[Definition~2.27]{Hackbusch} of $X+tZ$ is
\begin{align*}
&X+tZ =\\
&\begin{bmatrix} U_{:, 1:\ushort{r}-1} & U_{:, \ushort{r}+2:r+1} & U_{:, \ushort{r}:\ushort{r}+1}[w_+(\tau) \; -w_-(\tau)] \end{bmatrix}\\
&\,\diag(\sigma_1, \dots, \sigma_{\ushort{r}-1}, {\textstyle \frac{3}{4}} \tau \sigma_{\ushort{r}} I_{r-\ushort{r}}, \sigma_{\ushort{r}} \lambda_+(\tau), - \sigma_{\ushort{r}} \lambda_-(\tau))\\
&\begin{bmatrix} V_{:, 1:\ushort{r}-1} & V_{:, \ushort{r}+2:r+1} & V_{:, \ushort{r}:\ushort{r}+1}[w_+(\tau) \; w_-(\tau)] \end{bmatrix}^\tp
\end{align*}
if $\tau \in (\frac{12}{17}, 1)$, and
\begin{align*}
&X+tZ =\\
&\begin{bmatrix} U_{:, 1:\ushort{r}-1} & U_{:, \ushort{r}+2:r+1} & U_{:, \ushort{r}:\ushort{r}+1}[-w_-(\tau) \; w_+(\tau)] \end{bmatrix}\\
&\,\diag(\sigma_1, \dots, \sigma_{\ushort{r}-1}, {\textstyle \frac{3}{4}} \tau \sigma_{\ushort{r}} I_{r-\ushort{r}}, - \sigma_{\ushort{r}} \lambda_-(\tau), \sigma_{\ushort{r}} \lambda_+(\tau))\\
&\begin{bmatrix} V_{:, 1:\ushort{r}-1} & V_{:, \ushort{r}+2:r+1} & V_{:, \ushort{r}:\ushort{r}+1}[w_-(\tau) \; w_+(\tau)] \end{bmatrix}^\tp
\end{align*}
if $\tau \in (1, \frac{4}{3})$.
Therefore, by the Eckart--Young theorem,
\begin{align*}
&\proj{\R_{\le r}^{m \times n}}{X+tZ} =\\
&\begin{bmatrix} U_{:, 1:\ushort{r}-1} & U_{:, \ushort{r}+2:r+1} & U_{:, \ushort{r}:\ushort{r}+1} w_+(\tau)\end{bmatrix}\\
&\,\diag(\sigma_1, \dots, \sigma_{\ushort{r}-1}, {\textstyle \frac{3}{4}} \tau \sigma_{\ushort{r}} I_{r-\ushort{r}}, \sigma_{\ushort{r}} \lambda_+(\tau))\\
&\begin{bmatrix} V_{:, 1:\ushort{r}-1} & V_{:, \ushort{r}+2:r+1} & V_{:, \ushort{r}:\ushort{r}+1} w_+(\tau)\end{bmatrix}^\tp
\end{align*}
if $\tau \in (\frac{12}{17}, 1)$, and
\begin{align*}
&\proj{\R_{\le r}^{m \times n}}{X+tZ} =\\
&\begin{bmatrix} U_{:, 1:\ushort{r}-1} & U_{:, \ushort{r}+2:r+1} & - U_{:, \ushort{r}:\ushort{r}+1} w_-(\tau)\end{bmatrix}\\
&\,\diag(\sigma_1, \dots, \sigma_{\ushort{r}-1}, {\textstyle \frac{3}{4}} \tau \sigma_{\ushort{r}} I_{r-\ushort{r}}, - \sigma_{\ushort{r}} \lambda_-(\tau))\\
&\begin{bmatrix} V_{:, 1:\ushort{r}-1} & V_{:, \ushort{r}+2:r+1} & V_{:, \ushort{r}:\ushort{r}+1} w_-(\tau) \end{bmatrix}^\tp
\end{align*}
if $\tau \in (1, \frac{4}{3})$.
Since
\begin{align*}
&{U_{:, \ushort{r}:\ushort{r}+1}}^\tp
\left(\lim_{\substack{\tau \to 1 \\ <}} \proj{\R_{\le r}^{m \times n}}{X+tZ}\right)
V_{:, \ushort{r}:\ushort{r}+1}\\
&= \sigma_{\ushort{r}} \lambda_+(1) w_+(1) w_+(1)^\tp\\
&= \frac{\sigma_{\ushort{r}}}{4} \begin{bmatrix} 1 & 1 \\ 1 & 1 \end{bmatrix}\\
&\ne \frac{\sigma_{\ushort{r}}}{4} \begin{bmatrix} -1 & 1 \\ 1 & -1 \end{bmatrix}\\
&= \sigma_{\ushort{r}} \lambda_-(1) w_-(1) w_-(1)^\tp\\
&= {U_{:, \ushort{r}:\ushort{r}+1}}^\tp
\left(\lim_{\substack{\tau \to 1 \\ >}} \proj{\R_{\le r}^{m \times n}}{X+tZ}\right)
V_{:, \ushort{r}:\ushort{r}+1},
\end{align*}
every selection of the correspondence defined in~\eqref{eq:ProjectionTangentDeterminantalVariety} is discontinuous at $t_*$.
\hfill$\qed$
\end{pf}

\section{On a retracted line-search method for optimization on closed sets}
\label{sec:RetractedLineSearchClosedSet}
Given a continuously differentiable function $f : \mathcal{E} \to \R$, this section concerns the problem
\begin{equation}
\label{eq:MinDiffFunctionClosedSet}
\min_{x \in C} f(x)
\end{equation}
of minimizing $f$ on $C$. In general, without further assumptions on $C$ or $f$, problem \eqref{eq:MinDiffFunctionClosedSet} is intractable and algorithms are only expected to find a stationary point of this problem. A point $x \in C$ is said to be \emph{Bouligand stationary (B-stationary)} for~\eqref{eq:MinDiffFunctionClosedSet} if $\ip{\nabla f(x)}{v} \ge 0$ for all $v \in \tancone{C}{x}$ \cite[Definitions~6.1.1 and~1.1.3]{CuiPang}.

Line-search methods form an important class of algorithms for unconstrained optimization, i.e., optimization on a Euclidean vector space \cite[Chapter~3]{NocedalWright}. At every iteration, the next iterate is obtained by moving from the current iterate $x \in \mathcal{E}$ along a descent direction, i.e., a vector $v \in \mathcal{E}$ such that $\ip{\nabla f(x)}{v} < 0$.
Thanks to retractions, line-search methods have been extended from Euclidean vector spaces to Riemannian manifolds \cite[Chapter~4]{AbsilMahonySepulchre} and closed subsets of Euclidean vector spaces \cite[\S 2.4]{SchneiderUschmajew2015}. The resulting methods are called retracted line-search methods.

Algorithm~\ref{algo:BacktrackingRetractedLineSearch} summarizes the retracted line search performed at every iteration of the method proposed in \citet[Algorithm~1]{SchneiderUschmajew2015} for problem~\eqref{eq:MinDiffFunctionClosedSet}. It is based on the fact that, for every point $x \in C$ that is not B-stationary for~\eqref{eq:MinDiffFunctionClosedSet}, a descent direction $v \in \tancone{C}{x}$ exists.

\begin{algorithm}[H]
\caption{Backtracking retracted line search}
\label{algo:BacktrackingRetractedLineSearch}
\begin{algorithmic}[1]
\Require
$(\mathcal{E}, C, R, f, \beta, c)$ where $\mathcal{E}$ is a Euclidean vector space, $C$ is a nonempty closed subset of $\mathcal{E}$, $R$ is a retraction on $C$ (in the sense of Definition~\ref{defi:Retraction}), $f : \mathcal{E} \to \R$ is continuously differentiable, and $\beta, c \in (0, 1)$.
\Input
a point $x \in C$ that is not B-stationary for~\eqref{eq:MinDiffFunctionClosedSet}.
\Output
a point in $C$.

\State
Choose $v \in \tancone{C}{x}$ such that $\ip{\nabla f(x)}{v} < 0$;
\State
Choose $\alpha \in (0, \infty)$;
\label{algo:BacktrackingRetractedLineSearch:InitialStepSize}
\While
{$f(R(x, \alpha v)) > f(x) + c \, \alpha \, \ip{\nabla f(x)}{v}$}
\State
$\alpha \gets \alpha \beta$;
\EndWhile
\State
Return $R(x, \alpha v)$.
\end{algorithmic}
\end{algorithm}

Proposition~\ref{prop:BacktrackingRetractedLineSearchArmijo} states that the while loop in Algorithm~\ref{algo:BacktrackingRetractedLineSearch} terminates, thereby showing that the result in \citet[Proposition~2.8]{SchneiderUschmajew2015} holds for a retraction in the sense of Definition~\ref{defi:Retraction} which is less restrictive than that of \citet[Definition~2.4]{SchneiderUschmajew2015}. Specifically, the mere differentiability at $0$ of the function $[0, \infty) \to C : t \mapsto R(x, tv)$ for every $x \in C$ and $v \in \tancone{C}{x}$ ensures that the while loop terminates; the continuity of this function on $[0, \infty)$ is not necessary.

\begin{prop}
\label{prop:BacktrackingRetractedLineSearchArmijo}
Consider Algorithm~\ref{algo:BacktrackingRetractedLineSearch}.
There exists $\alpha_* \in (0, \infty)$ such that, for all $t \in [0, \alpha_*]$, the Armijo condition
\begin{equation*}
f(R(x, tv)) \le f(x) + c \, t \, \ip{f(x)}{v}
\end{equation*}
is satisfied.
Thus, the while loop terminates after at most $\max\{0, \lceil\ln(\alpha_*/\alpha)/\ln(\beta)\rceil\}$ iterations with a step size $\ushort{\alpha}$ at least $\min\{\alpha, \alpha_* \beta\}$, where $\alpha$ is the initial step size chosen in line~\ref{algo:BacktrackingRetractedLineSearch:InitialStepSize}.
\end{prop}

\begin{pf}
Since $R$ is a retraction on $C$, the function
\begin{equation*}
\hat{R} : [0,\infty) \to C : t \mapsto R(x, tv)
\end{equation*}
satisfies $\hat{R}(0) = x$ and $\hat{R}'(0) = v$. Thus, by the chain rule,
\begin{equation*}
(f \circ \hat{R})'(0) = \ip{\nabla f(x)}{v}.
\end{equation*}
Let $\varepsilon \coloneq - (1-c) \ip{\nabla f(x)}{v}$. There exists $\alpha_* \in (0, \infty)$ such that, for all $t \in (0, \alpha_*]$,
\begin{equation*}
\left|\frac{(f \circ \hat{R})(t)-(f \circ \hat{R})(0)}{t}-(f \circ \hat{R})'(0)\right| \le \varepsilon,
\end{equation*}
which implies
\begin{equation*}
\frac{(f \circ \hat{R})(t)-(f \circ \hat{R})(0)}{t} \le c \, \ip{\nabla f(x)}{v},
\end{equation*}
i.e., the Armijo condition.
Therefore, either the initial step size $\alpha$ chosen in line~\ref{algo:BacktrackingRetractedLineSearch:InitialStepSize} satisfies the Armijo condition or the while loop ends after iteration $i \in \N \setminus \{0\}$ with $\ushort{\alpha} = \alpha\beta^i$ such that $\ushort{\alpha}/\beta > \alpha_*$. In the second case, $i < 1+\ln(\alpha_*/\alpha)/\ln(\beta)$ and thus $i \le \lceil\ln(\alpha_*/\alpha)/\ln(\beta)\rceil$.
\hfill$\qed$
\end{pf}

\section{Conclusion}
\label{sec:Conclusion}
This paper brings three contributions to the study of retractions on closed subsets of Euclidean vector spaces. First, projective retractions are proven to exist on geometrically derivable sets (Proposition~\ref{prop:ProjectiveRetractionGeometricDerivability}), a class of sets that contains algebraic varieties. Second, it is shown that no projective retraction on the determinantal variety $\R_{\le r}^{m \times n}$ is a retraction in the sense of \citet[Definition~2.4]{SchneiderUschmajew2015} (Proposition~\ref{prop:ProjectiveRetractionDeterminantalVariety}). Third, Algorithm~\ref{algo:BacktrackingRetractedLineSearch} is proven to terminate (Proposition~\ref{prop:BacktrackingRetractedLineSearchArmijo}), thereby ensuring an Armijo condition, a result already stated in \citet[Proposition~2.8]{SchneiderUschmajew2015} for the stronger notion of retraction proposed in \citet[Definition~2.4]{SchneiderUschmajew2015}. The convergence analysis of several algorithms aiming at solving problem~\eqref{eq:MinDiffFunctionClosedSet} with $C$ a low-rank variety would benefit from these results. Examples of such algorithms include the method proposed in \citet[Algorithm~3]{SchneiderUschmajew2015} for $C = \R_{\le r}^{m \times n}$, known as $\ppgd$, and the GRAP method proposed in \citet[Algorithm~2]{GaoPengYuan} for $C$ the algebraic variety of tensors of upper-bounded multilinear rank.
These three contributions suggest that, in the context of constrained optimization, the definition of retraction proposed in \citet[\S 3.1.2]{HosseiniUschmajew2019} should be preferred to that proposed in \citet[Definition~2.4]{SchneiderUschmajew2015} as it inherits the main property while being less restrictive.

\section*{Acknowledgements}
The author thanks P.-A. Absil and two anonymous referees for several helpful comments that improved the quality of the paper.

\bibliography{golikier_bib}
\end{document}